\documentclass[preprint]{revtex4-1}
\usepackage{amsmath,amssymb}
\usepackage{amsthm}
\usepackage{graphicx}
\usepackage{hyperref}
\usepackage{xcolor}

\newcommand{\bd}[1]{\boldsymbol{#1}}
\newcommand{\wt}[1]{\widetilde{#1}}

\newcommand{\ie}{{\it i.e.}}
\newcommand{\eg}{{\it e.g.}}

\renewcommand{\i}{\ensuremath{\text{\normalfont I}}}
\newcommand{\ii}{\ensuremath{\text{\normalfont I\!I}}}
\newcommand{\ici}{\ensuremath{\text{\normalfont I,I}}}

\newcommand{\iici}{\ensuremath{\text{\normalfont I\!I,I}}}
\newcommand{\iicii}{\ensuremath{\text{\normalfont I\!I,I\!I}}}

\begin{document}

\title{Absorbing Boundary Conditions for Time-dependent Schr\"{o}dinger equations: A Density-matrix Formulation}

\author{ Xiantao Li}
\email{xli@math.psu.edu}
\affiliation{Department of Mathematics, the Pennsylvania State University, University Park, PA 16802-6400, USA.}%

\date{\today}

\begin{abstract} 
This paper presents some absorbing boundary conditions (ABC) for simulations based on the time-dependent density-functional theory (TDDFT). 
The boundary conditions are expressed in terms of the elements of the density-matrix, and it is  derived from the full model over a much larger domain.
To make the implementation much more efficient, several approximations for the convolution integral  will be constructed with guaranteed stability. These approximations lead
to modified density-matrix equations at the boundary. The effectiveness is examined via numerical tests. 
\end{abstract} 

\maketitle

\section{Introduction}

Time-dependent density-functional theory  (TDDFT)  \cite{marques2004time,martin2004,ullrich2011time} has recently become an extremely useful 
computational tool to study electronic and optical properties of materials and bio-molecules. The model  describes the many-body problem by using a noninteracting system 
that  reproduces the electronic density of the full, interacting system via an exchange-correlation functional \cite{runge1984density}. There have been tremendous progress in the development
of  efficient algorithms for the implementation of TDDFT\cite{ullrich2011time}. 
One  computational issue emerges, when electrons excited to the unbound states are emitted outside the system. The usual periodic or
Dirichlet boundary conditions will simply introduce the electrons back into the system, which will subsequently interfere with the computation. Simulating such processes  
accurately demands  
an efficient absorbing boundary condition (ABC) that  suppresses the artificial  boundary reflections \cite{yabana2006real}. 

In essence, ABCs are domain reduction methods: The exact model is defined in a relatively large domain $\Omega$, and a direct simulation over this domain is intractable due to the overwhelmingly large number of degrees of freedom.  Instead, one aims to reduce the problem to a sub-domain, say, $\Omega_\i$,  $\Omega_\i \subset \Omega$, by deriving an appropriate boundary condition at the boundary of $\Omega_\i$. As a result, the computation can be performed over a much smaller domain, at a reduced cost. The accuracy of the ABCs is often  manifested in the magnitude of wave reflections at the boundary. This has prompted researchers to formulate the ABCs by minimizing the reflection coefficients. ABCs have a wide variety of applications, ranging from acoustic wave equations \cite{Alpert2002,Engquist1977,Hagstrom1999}, Maxwell equation \cite{berenger1994}, Helmholtz equation \cite{Oberai1998}, molecular dynamics \cite{Cai2000,Karpov2005,Li08b,Li2006,E2001,Wu2018}, etc. 

In the context of time-dependent Schr\"{o}dinger equations, there are two important cases where boundary conditions can be useful. The first case is when the surrounding medium is a vacuum. Namely, the electron density is  initially ($t=0$) completely concentrated in the domain $\Omega_\i$. By further assuming  that the external potential is also confined in the same domain, one can derive an ABC for the wave function, typically expressed as a time convolution at the boundary. Such a problem has been extensively studied in the applied math community \cite{antoine2008review,antoine2004,arnold1998,arnold2003,fevens1999,han2004,han2005finite,Jiang2004,Jiang2008,xu2006absorbing,zheng2007}. Meanwhile, there are several techniques that are specifically designed for TDDFT models, including the mask function\cite{krause1992calculation},  imaginary potentials \cite{child1991analysis,muga2004complex}, which has resemblance to the perfectly matched layer method \cite{berenger1994,Mur1981}, as well as the Green's function approach \cite{yabana2006real}.

Another important case is when the system under consideration is part of a bigger quantum-mechanical system, where there are also electrons outside the computational domain  $\Omega_\i$. Consequently, the goal of the ABCs is to start with a large domain with a large number of electrons (and many equations), and derive an effective  model in a much smaller domain with much fewer electrons (and fewer equations).  One example, is the interaction between a crystalline solids and an ion beam \cite{hatcher2008dynamical}, 
where the quantum description can be confined to a localized region. The main challenge in this case is that the wave functions are often extended, due to the Bloch theorem, and a premature truncation would lead to large error. 

In this paper, we will consider the second class of problems. In particular, we  formulate the ABCs in terms of the density-matrix. 
We will show that the ABC can be expressed in terms of the components of the density-matrix at the boundary of the sub-domain $\Omega_\i$, in the form of a convolution. Furthermore, the matrix function in the convolution can be approximated in such a way that the integral does not need to be repeatedly computed. We also prove that the resulting boundary condition is always stable.

The motivation behind a density-matrix formulation is { threefold}. First, the density-matrix is represented directly at the grid points and it eliminates the explicit dependence on the type/number of   local orbitals, which makes the derivation easier. Second, the resulting ABCs can be directly applied to density-matrix implementation  of the TDDFT model  \cite{furche2001density,xie2012time,yam2003localized}, which can be made quite efficient by making use of the nearsightedness of the density-matrix. Finally, 
 our approximation schemes might be useful for simulating non-Markovian effects arising from open quantum systems \cite{groblacher2015observation,madsen2011observation,mi2016strong,strathearn2018efficient}.

The remaining part of the paper is organized as follows.  Sec. \ref{th} presents the setup of the problem and the derivation of the boundary condition in terms of the elements of the density-matrix at the boundary. In Sec. \ref{ap}, several approximation schemes will be constructed based on the Laplace transform of the matrix function.  To examine the effectiveness of the approximations, we present results from some numerical experiments in Sec. \ref{num}. The test problem has been motivated by the simulations \cite{de2012wave,gruber2016ultrafast}  for a graphene sheet under localized external field.

\section{The basic theory}\label{th}

In TDDFT, the underlying description is the time-dependent Schr\"{o}dinger equations,
\begin{equation}\label{eq: schr}
 i\partial_t \psi_\ell = \hat{H} \psi_\ell, \quad \ell=1,2, \cdots, N_e.
\end{equation}
Here $N_e$ is the total number of electrons. $\hat{H}$ is the Hamiltonian operator, consisting of the kinetic energy, an effective potential and an external potential,
\begin{equation}
 \hat{H} = \hat{T} + \hat{V}_{\text{eff}} + \hat{V}_{\text{ext}}(\bd r, t).
\end{equation} 

We will denote $\hat{H}^{(0)} = \hat{T} + \hat{V}_{\text{eff}}$ as an unperturbed Hamiltonian. Typically, the initial condition of the wave functions are prepared as the ground states of $\hat{H}^{(0)}$ in the absence of the external potential. Then by activating the external potential,  equations \eqref{eq: schr} can be solved to study the properties in response to $V_{\text{ext}}$, \eg, the polarizability and absorption. In practice, the computational cost grows rapidly as the number of electrons increases.  This paper considers a scenario where such computation can be reduced by using appropriate boundary conditions. 

\subsection{The derivation of the absorbing boundary condition}

Our starting point is the dynamics in terms of the density-matrix. More specifically, we define the density-matrix operator,
\begin{equation}\label{eq: dm}
  \hat{\rho}(\bd r, \bd r', t)= \sum_\ell n_\ell \psi_\ell(\bd r, t) \psi_\ell(\bd r', t)^*,
\end{equation}
with $n_\ell$'s being the occupation numbers.  

The density-matrix satisfies the Liouville-von Neumann equation
\begin{equation}\label{eq: rho}
i \partial_t \hat{\rho} =\big[\hat{H}, \hat{\rho} \big].
\end{equation} 
Here $\rho(\bd r, \bd r', t)$ is the density-matrix with $\bd r, \bd r' \in \Omega$; $\Omega$ is the physical domain for the entire system.   


To formulate the absorbing boundary condition (ABC), we first assume that equation \eqref{eq: rho} has been appropriately discretized, so that $\rho$ is a matrix defined at certain grid points. As a result, one obtains a matrix-valued (finite-dimensional) system, and hence we will drop the $\hat{\cdot}$  notation from now on. We now assume that the external potential  $V_{\text{ext}}$ is only non-zero in a sub-domain, denoted here by $\Omega_I$; The surrounding region is denoted by $\Omega_{\ii}$; $\Omega_I \cup \Omega_{\ii}=\Omega.$ We also set up the problem by assuming that initially the system is at ground states, for which  the density-matrix is denoted by $\rho^0.$ In particular, we have $\big[H^{(0)}, \rho^0\big]=0.$

In general, equation \eqref{eq: rho} is nonlinear, since $H$ may depend on $\rho$. For instance, in the TDDFT model, the effective potential $V_{\text{eff}}$ in the Hamiltonian depends on the electron density. Linearization can be made around the ground states, which has been the starting point of the linear response theory \cite{casida1995time,muta2002solving,yabana2006real}. {For instance,} the Hamiltonian can be linearized to,
\begin{equation}
  [H, \rho] \approx   [H^{(0)}, \delta \rho]  +  \big[ \int \frac{\delta V_{\text{eff}}[n_0(\bd r,t)]}{\delta n(\bd r', t)} \delta n(\bd r', t) d\bd r',   \rho^0 \big]
\end{equation}
Here we have used $n(\bd r, t)$ for the electron density,  $n_0(\bd r)$ as the ground-state electron density,  and $\delta\!n(\bd r, t)= n(\bd r, t) -n_0(\bd r) $ as the 
corresponding perturbation perturbation.  Typically in the formulation of absorbing boundary conditions, this is only required in the exterior domain, and the nonlinearity can be retained 
in $\Omega_I$ \cite{Wu2018}.
For simplicity, we will omit the second term in this paper and use $H^{(0)}$ in the density-matrix equation.

\begin{figure}[htbp]
\begin{center}
\includegraphics[scale=0.45]{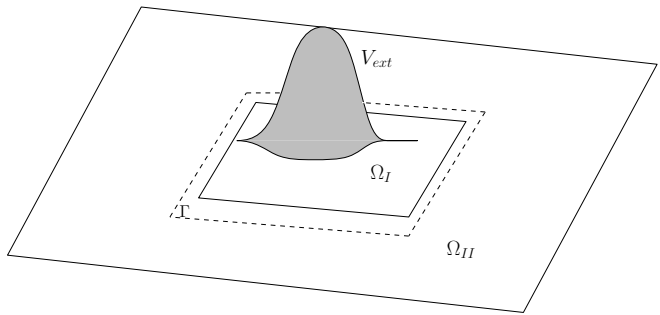}
\caption{An illustration of the problem setup for the derivation of absorbing boundary conditions. The goal is to reduce the problem to a subdomain ($\Omega_\i$), where the external potential is confined, by eliminating the surrounding electronic degrees of freedom (in $\Omega_\ii).$ The resulting boundary conditions will be written in terms of the density-matrix at the artificial boundary ($\Gamma$). } 
\label{default}
\end{center}
\end{figure}

In accordance with the partition of the domain, we will partition the density-matrix and the Hamiltonian operator. For example, we denote $\rho_{\ici}$ as the restriction of $\rho$ to the domain $\Omega_\i$, \ie, $\bd r \in \Omega_\i$ and  $\bd r' \in \Omega_\i$. This separates equation \eqref{eq: rho} into three equations,  given by,
\begin{equation}\label{eq: rho'}
\begin{aligned}
 i \partial_t \rho_{\i,\i} =& H_{\i,\i}\rho_{\i,\i} + H_{\i,\ii}\rho_{\ii,\i} - \rho_{\i,\i}H_{\i,\i}- \rho_{\i,\ii}H_{\ii,\i},\\
 i \partial_t \rho_{\i,\ii} =& H_{\i,\i}\rho_{\i,\ii} + H_{\i,\ii}\rho_{\ii,\ii} - \rho_{\i,\i}H_{\i,\ii}- \rho_{\i,\ii}H_{\ii,\ii},\\
  i \partial_t \rho_{\ii,\ii} =& H_{\ii,\i}\rho_{\i,\ii} + H_{\ii,\ii}\rho_{\ii,\ii} - \rho_{\ii,\i}H_{\i,\ii}- \rho_{\ii,\ii}H_{\ii,\ii}. 
\end{aligned}
\end{equation}
The block entry $\rho_{\ii,\i}$ is given by $\rho_{\i,\ii}^*$ due to its Hermitian property. The partition of the domain is illustrated in Figure \ref{default}.

%

We recall that initially, the density-matrix in the exterior is at equilibrium,
\begin{equation}
 \rho_{\ii,\ii}(0)= \rho_{\ii,\ii}^0.
\end{equation}
Our derivation will concentrate on the perturbation, 
\begin{equation}
\delta\!\rho_{\ii,\ii}(t)= \rho_{\ii,\ii}(t)- \rho_{\ii,\ii}^0,
\end{equation}
 as a response to the dynamics in the interior $\Omega_\i$. 

By applying the variation-of-constant formula, we find,
\begin{equation}\label{eq: rho22}
\delta\!\rho_{\ii,\ii}(t) =-i\int_0^t e^{H_{\ii,\ii}(t-t')/i} \Big(H_{\ii,\i}\delta\!\rho_{\ii,\i}(t')- \delta\!\rho_{\ii,\i}(t')H_{\i,\ii}\Big)e^{-H_{\ii,\ii}(t-t')/i}dt'
\end{equation}

With this formula, we now turn to the second equation in \eqref{eq: rho'}, which now becomes 
\begin{equation}\label{eq: rho12}
 i \partial_t \rho_{\i,\ii} = H_{\i,\i}\rho_{\i,\ii} + H_{\i,\ii}\rho_{\ii,\ii}^0 + H_{\i,\ii}\delta\underline{\underline{\!\rho_{\ii,\ii}} }- \rho_{\i,\i}H_{\i,\ii}- \rho_{\i,\ii}H_{\ii,\ii},\\
\end{equation}

In principle, equation \eqref{eq: rho22} can be substituted into equation \eqref{eq: rho12} so that the large-dimensional system in $\Omega_{\ii}$ is eliminated. However, the number of elements in $\rho_{\i,\ii}$ is still large. To achieve further reduction, we assume that the Hamiltonian matrix is localized, in that when two grid points are far apart, the corresponding entry is zero. For example, 
 the kinetic energy can be approximated by a finite-difference method with a 3-point or  5-point stencil in each direction.     
To take advantage of the locality of $H$, we define a boundary region $\Gamma \subset \Omega_{\ii}$ 
commensurate with
 the  cut-off distance of the Hamiltonian matrix.  We further define a restriction operator from $\Omega_{\ii}$ to $\Gamma$, denoted by $R$. 
 As a result, we have
 \begin{equation}
 H_{\Gamma,\ii} =R H_{\iicii}, \; H_{\Gamma,\i}=R H_{\iici}.
\end{equation}

Similarly, we define part of the density matrices near the boundary as follows,
\begin{equation}
  \rho_{\Gamma,\ii} = R \rho_{\iicii},  \;  \rho_{\Gamma,\i}= R \rho_{\ii,\i}.
  \end{equation}

Two  conditions can be deduced from this matrix: $RR^T=I$, with $I$ being the identity operator on $\Omega_\Gamma,$ and, $R^T R$ is an orthogonal projection operator on $\Omega_\ii$.  We will denote the dimension of the subdomains $\Omega_\i$, $\Omega_\ii$ and $\Gamma$, by $n_\i$, $n_\ii$ and $n_\Gamma$, respectively. $n_\i, n_\Gamma \ll n_\ii$.

With a multiplication by $R^T$ from the right, equation \eqref{eq: rho12} can be written as,
\begin{equation}\label{eq: rho12'}
 i \partial_t \rho_{\i,\Gamma} = H_{\i,\i}\rho_{\i,\Gamma} + H_{\i,\Gamma}\rho_{\Gamma,\Gamma}^0 + H_{\i,\Gamma}\underline{\underline{\delta\!\rho_{\Gamma,\Gamma} }}- \rho_{\i,\i}H_{\i,\Gamma}- \rho_{\i,\Gamma}H_{\Gamma,\Gamma}.
\end{equation}
Here we have used the fact that $H_{\i,\ii} \rho_{\ii,\ii}=H_{\i,\Gamma} R \rho_{\ii,\ii}=  H_{\i,\Gamma}  \rho_{\Gamma,\ii}$.  By keeping just the entries in $\rho_{\i,\Gamma}$, many columns in the 
matrix $\rho_{\i,\ii}$ in \eqref{eq: rho12} have been removed.

Using the same observation, one can simplify the first equation in \eqref{eq: rho'} as follows 
\begin{equation}\label{eq: rho11'}
 i \partial_t \rho_{\i,\i} = H_{\i,\i}\rho_{\i,\i} + H_{\i,\Gamma}\rho_{\Gamma,\i} - \rho_{\i,\i}H_{\i,\i}- \rho_{\i,\Gamma}H_{\Gamma,\i}.
\end{equation}

Equations \eqref{eq: rho11'} and \eqref{eq: rho12'} can be readily solved, provided that $\delta\!\rho_{\Gamma,\Gamma}(\cdot,\cdot,t)$ is known. Therefore, the  ABCs will be expressed in terms of   $\delta\!\rho_{\Gamma,\Gamma}(\cdot,\cdot,t)$ (underlined in the equation \eqref{eq: rho12'}). By applying the restriction operator to \eqref{eq: rho22}, we obtain,
\begin{equation}
\delta \rho_{\Gamma,\Gamma}(t) = -i\int_0^t R e^{H_{\ii,\ii}(t-t')/i}  \Big(R^T H_{\Gamma,\i}\delta\!\rho_{\i,\ii}(t')- \delta \rho_{\ii,\i}(t')  H_{\i,\Gamma}R\Big)  e^{-H_{\ii,\ii}(t-t')/i} R^T dt'.
\end{equation}

We make the truncation that $\delta\!\rho_{\i,\ii}\approx \delta\!\rho_{\i,\Gamma} R^T.$ This is motivated by the nearsightedness property of the density-matrix. This simplifies the above equation to a  closed-form formula,
\begin{equation}\label{eq: delta-rho-gg}
\delta\! \rho_{\Gamma,\Gamma}(t) = -i\int_0^t  Y(t-t') \Big(H_{\Gamma,\i}\delta\!\rho_{\i,\Gamma}(t')- \delta\!\rho_{\Gamma,\i}(t')  H_{\i,\Gamma}\Big)  Y(t-t')^* dt',
\end{equation}
where the matrix function $Y$ is given by,
\begin{equation}
 Y(t)= R e^{-iH_{\ii,\ii}t} R^T. 
\end{equation}

 At this point, our ABC is expressed as a convolution integral with kernel function given by $Y(t).$ The function $Y(t)$ has dimension $n_\Gamma \times n_\Gamma$.
However,  direct calculations would still require the computation of a matrix exponential involving a matrix with dimension
 $n_\ii \times n_\ii$, which would be expensive to compute directly. Next we discuss methods to simplify this calculation, and make the overall implementation more efficient.

\section{Approximations of the absorbing boundary condition }\label{ap}

The crudest approximation is to neglect the influence of $\delta\! \rho_{\Gamma,\Gamma}(t)$ in \eqref{eq: delta-rho-gg} entirely by simply setting it to zero.  As a result, we solve,
\begin{equation}\label{eq: rho12-d}
 i \partial_t \rho_{\i,\Gamma} = H_{\i,\i}\rho_{\i,\Gamma} + H_{\i,\Gamma}\rho_{\Gamma,\Gamma}^0  - \rho_{\i,\i}H_{\i,\Gamma}- \rho_{\i,\Gamma}H_{\Gamma,\Gamma}.
\end{equation}
In this case, one takes into account the surrounding environment via the term \( H_{\i,\Gamma}\rho_{\Gamma,\Gamma}^0 \). However, there is no feedback from the region 
$\Omega_{\ii}.$ Therefore, this approximation will not be considered here. Instead, we construct more accurate approximations based on  the Laplace transform of $Y(t)$, given by,
\begin{equation}
  \wt{Y} (s)= \int_0^{+\infty} Y(t) e^{-st}dt = R\big[sI + iH_{\ii,\ii}\big]^{-1}R^T, \quad \text{Re}(s)>0.
\end{equation}

Next, we present some approximation methods. 

\subsection{A first-order approximation}
As a first-order approximation, we approximate  $Y(t)$ by a delta function,
\begin{equation}
  Y(t)  \approx Y_0 \delta(t). 
\end{equation}
Here we choose 
\begin{equation}
 Y_0 = \wt{Y}(s_0),
\end{equation}
which only needs to be computed once.

As a result, the boundary condition \eqref{eq: delta-rho-gg} becomes,
\begin{equation}
\delta \!\rho_{\Gamma,\Gamma}(t)  = -i Y_0 \Big(H_{\Gamma,\i}\delta\!\rho_{\i,\Gamma}(t)- \delta\!\rho_{\Gamma,\i}(t)  H_{\i,\Gamma}\Big) Y^*_0,
\end{equation}
which can be substituted into \eqref{eq: rho12'} to complete the model. This approximation will be referred to as the {\bf first-order} ABC.

\subsection{Second-order approximations}

The previous method approximates the Laplace transform $\wt{Y}(s)$ by a constant matrix. To derive a more accurate approximation,  we first make the observation that 
\begin{equation}
  Y(0)=I. 
\end{equation}
Therefore, we will seek an  approximation as a matrix exponential,
\begin{equation}\label{eq: zt}
 Z(t)= e^{-i\bar{H}_{\Gamma,\Gamma} t},
\end{equation}
with $\bar{H}_{\Gamma,\Gamma}$ being a $n_\Gamma \times n_\Gamma$ matrix that will be determined as follows.

The choice is constructed to first satisfy the obvious condition $Z(0)=I.$ To determine the matrix $\bar{H}_{\Gamma,\Gamma}$, we set an interpolation condition, 
\begin{equation}\label{eq: interp0}
 \wt{Z}(s_0)= \wt{Y}(s_0). 
\end{equation}
This gives the expression for $\bar{H}_{\Gamma,\Gamma}$, 
\begin{equation}\label{eq: Hgg}
 \bar{H}_{\Gamma,\Gamma} =(-i)\big( \wt{Y}(s_0)^{-1} - s_0 I \big).
\end{equation}

In light of equations \eqref{eq: zt} and \eqref{eq: delta-rho-gg}, $\rho_{\Gamma,\Gamma}$ follows the effective dynamics with no history dependence,
\begin{equation}\label{eq: rhoGG}
i\partial_t\delta\! \rho_{\Gamma,\Gamma}= [\bar{H}_{\Gamma,\Gamma}, \delta\! \rho_{\Gamma,\Gamma}] 
 + H_{\Gamma,\i}\delta\!\rho_{\i,\Gamma}- \delta\!\rho_{\Gamma,\i}  H_{\i,\Gamma},
\end{equation}
which provides an alternative for the numerical implementation.

Notice that the matrix  $\bar{H}_{\Gamma,\Gamma}$ with dimension $n_\Gamma \times n_\Gamma$  only needs to be computed once. Then equation \eqref{eq: rhoGG} can be solved at each
time step to provide the necessary boundary conditions for the dynamics in the interior  ( Equations \eqref{eq: rho11'} and \eqref{eq: rho12'}). Another observation is that
the combined dynamics  \eqref{eq: rho11'}, \eqref{eq: rho12'} and \eqref{eq: rhoGG} can also be realized by solving 
the time-dependent Schr\"{o}dinger equations,
\begin{equation}
  i \partial_t
  \left[\begin{array}{c}\psi_I \\ \psi_\Gamma \end{array}\right] = \left[\begin{array}{cc}H_{I,I} & H_{I,\Gamma} \\ H_{\Gamma,I} & \bar{H}_{\Gamma,\Gamma} \end{array}\right]\left[\begin{array}{c}\psi_I \\ \psi_\Gamma \end{array}\right] + V_{ext} \left[\begin{array}{c}\psi_I \\ \psi_\Gamma \end{array}\right] .
\end{equation} 

We will referred the method   \eqref{eq: rhoGG} as the {\bf second-order IIa} ABC. Another choice is neglect the initial condition $Y(0)=I,$ and replace it by an extra interpolation condition,
\begin{equation}
  \wt{Z}'(s_0)= \wt{Y}'(s_0). 
\end{equation}

This is a Pad\`{e} approximation at $s=s_0,$ as motivated by the ABCs for the wave equation \cite{Engquist1977}. In the time domain, the approximation can be expressed as,
\begin{equation}
   Z(t)= e^{-i\bar{H}_{\Gamma,\Gamma} t}A,
\end{equation}
with the two matrices given by,
\begin{equation}
  \bar{H}_{\Gamma,\Gamma} = i \big(s_0I + \wt{Y}(s_0)\wt{Y}'(s_0)^{-1}\big), \quad A= \big(s_0I + i  \bar{H}_{\Gamma,\Gamma} \big) \wt{Y}(s_0).
\end{equation}
The effective dynamics  for $\rho_{\Gamma,\Gamma}$ is slightly different,
\begin{equation}\label{eq: rhoGG'}
i\partial_t\delta\! \rho_{\Gamma,\Gamma}= [\bar{H}_{\Gamma,\Gamma}, \delta\! \rho_{\Gamma,\Gamma}] 
 + A \big( H_{\Gamma,\i}\delta\!\rho_{\i,\Gamma}- \delta\!\rho_{\Gamma,\i}  H_{\i,\Gamma} \big)A^*.
\end{equation}
This method  \eqref{eq: rhoGG'} will be called  the {\bf second-order IIb} ABC.

\subsection{A Stability Analysis}

The equation \eqref{eq: rhoGG} might be seen as  a modification of the dynamics at the boundary $\Omega_\Gamma$, with the effective Hamiltonian
given by $\bar{H}_{\Gamma,\Gamma}$. First, we observe that  $\bar{H}_{\Gamma,\Gamma}$ is a symmetric matrix.  The dynamics  \eqref{eq: rhoGG} that describes the ABC is stable if the imaginary part of  $\bar{H}_{\Gamma,\Gamma}$ is semi negative-definite.  

We now show that,  when $s_0>0, $    $\bar{H}_{\Gamma,\Gamma}$ from the interpolation \eqref{eq: interp0} always satisfies the stability condition. To  examine the matrix, we first separate the real and imaginary parts of $\wt{Y}(s_0)$, 
\[ \wt{Y}(s_0) =s_0 R\big(s_0^2I + H_{\ii,\ii}^2 \big)^{-1}R^T - i R\big(s_0^2I + H_{\ii,\ii}^2 \big)^{-1}H_{\ii,\ii} R^T.\]
Therefore, the real part of $\wt{Y}(s_0)$ is a positive definite matrix when $s_0>0.$  We now turn to  $\bar{H}_{\Gamma,\Gamma}$ in  \eqref{eq: Hgg}  by writing,
\[  \wt{Y}(s_0)^{-1} -s_0 I =   \wt{Y}(s_0)^{-1}\Big( \wt{Y}(s_0)^* - s_0 \wt{Y}(s_0) \wt{Y}(s_0)^* \Big)  \wt{Y}(s_0)^{-*}.\]

It is enough to examine the real part of the matrix in the middle of  the term on the right hand side,
\[\begin{aligned}
\text{Re}& \Big(\wt{Y}(s_0)^* -s_0 \wt{Y}(s_0) \wt{Y}(s_0)^* \Big) \\ 
 &= R\big(s_0^2I + H_{\ii,\ii}^2 \big)^{-1}\Big(s_0^3I + s_0 H_{\ii,\ii}^2 -s_0^3 R^TR - s_0 H_{\ii,\ii} R^T R H_{\ii,\ii}\Big)\big(s_0^2I + H_{\ii,\ii}^2 \big)^{-1}R^T.
\end{aligned} \]
Using the fact that $R^TR$ is an orthogonal projection matrix, we conclude that the real part of this matrix is positive definite.  As a result,  the imaginary part of $\bar{H}_{\Gamma,\Gamma}$ a  semi negative-definite  matrix,  since
\[ \text{Im} \bar{H}_{\Gamma,\Gamma} = - \text{Re} \big(\wt{Y}(s_0)^{-1} -s_0 I\big) \leq 0.\]
Hence the stability requirement is fulfilled. 

It is possible to choose a complex value for $s_0$. In this case the stability will still follow as long as the real part of $s_0$ is positive, since the imaginary part can be combined with $H_{\ii,\ii}$ in \eqref{eq: Hgg}. For example, one can choose $\text{Im}s_0= \omega,$ which corresponds to the frequency of the external potential. 

Finally, when $s_0=0,$ $\bar{H}_{\Gamma,\Gamma}$ is real symmetric. So the stability condition is also satisfied. But in this case, the total number of electrons will be conserved, which would not offer the desired absorbing property. 

\section{Numerical results}\label{num}

As a numerical test, we consider a single graphene sheet, consisting of 160 atoms. This forms the entire system. The Hamiltonian matrix $H$ is generated in OCTOPUS \cite{marques2003octopus}, which implements a real-space method for the TDDFT model. Atomic units will be used, including the Bohr radius as the length unit and the time unit is given by $2.4188843 \time 10^{-17}$ second.  For the discretization of the Hamiltonian, we choose the grid size to be  $\Delta x = 0.3677$.  The system is set up so that it is periodic in the first two dimensions, with dimension $28.6842 \times 33.0971 \times 7.354.$  

The initial density-matrix is computed based on the eigenfunctions of $H$. We consider an external potential given by,
\begin{equation}
  V_{\text{ext}}(\bd r, t) = e^{-c \bd r^2} e^{-\gamma (t-t_0)^2} \sin \omega t.
\end{equation}
Here the parameters are chosen as follows. $c=0.35$, $\gamma =3,$ $t_0=1.5$, and $\omega=6$, in atomic units. The first term  focuses the potential at the center of the domain. 
The remaining two terms provide a pulse, centering the frequency around $\omega$. 

\begin{figure}[thp]
\begin{center}
\includegraphics[scale=0.25]{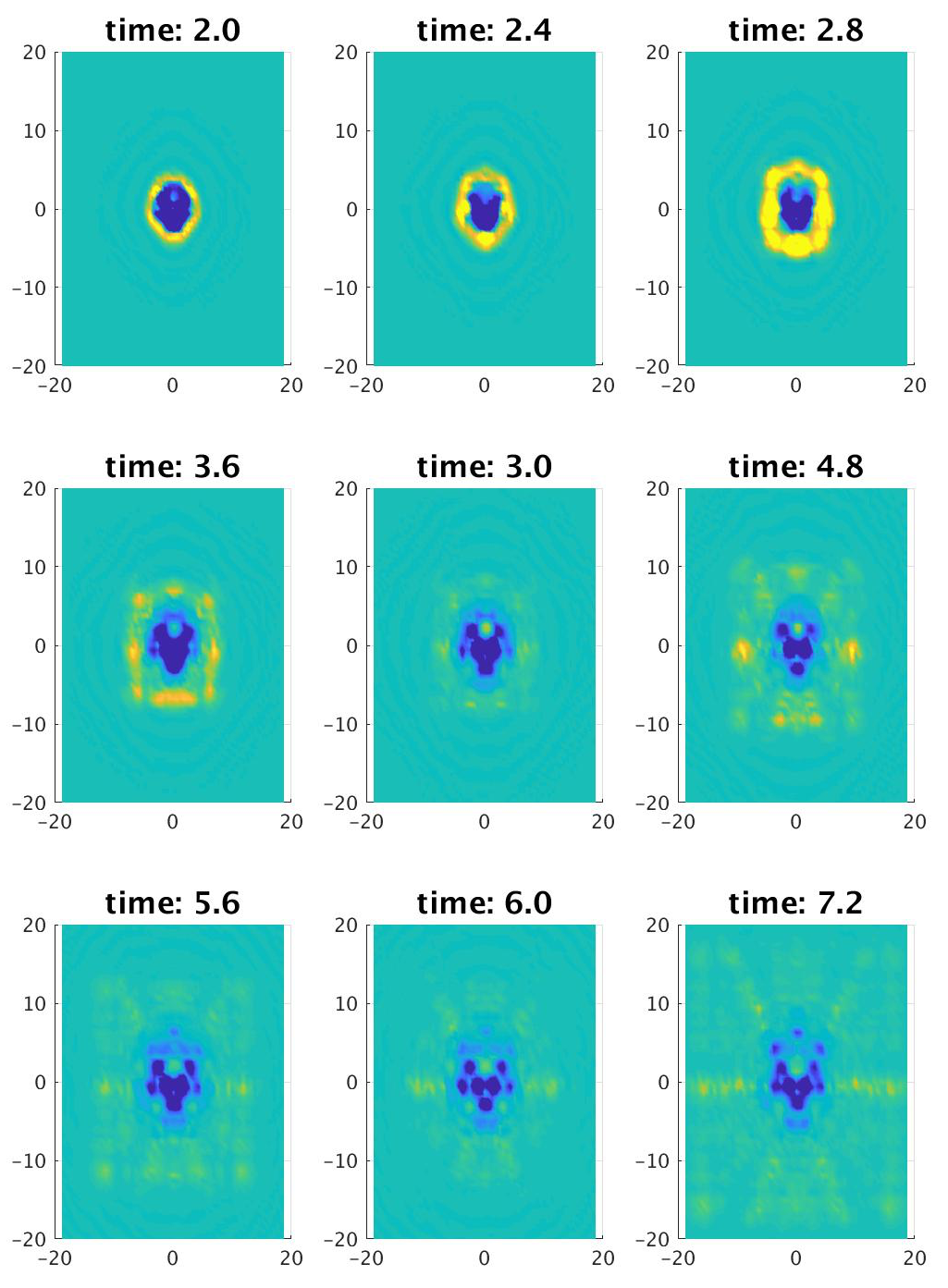}
\end{center}
\caption{The change of the electron density, $\delta n(\bd r, t),$ from the full simulation. }\label{eq: full}
\end{figure}

We first solve the full dynamics \eqref{eq: rho}.  Figure \ref{eq: full} shows snapshots for the change of the electron density, $ \delta n(\bd r, t),$ during the full simulation. We observe that wave front  of hexagon shape is generated at the center, which might be due to the lattice structure, and it begins to propagate toward the exterior. The wave front first  evolves into a
square shape, and then an octagon shape, followed by a square shape again. A dispersed pattern of waves is observed  later one. At around $t=7.2$, the waves have arrived at the outer boundaries.

\begin{figure}[tp]
\begin{center}
\includegraphics[scale=0.32]{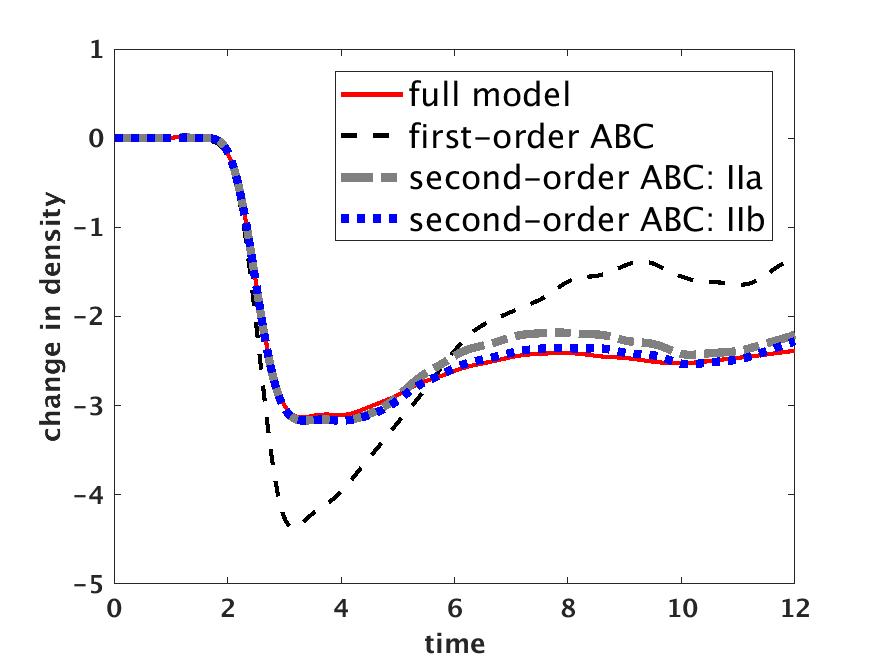}
\end{center}
\caption{The time history of the change of the total electrons in the subdomain $\Omega_{\i}=[-4,\;\;4] \times [-4,\;\;4]$. }\label{eq: ech}
\end{figure}

In order to test the ABCs, we choose the subdomain $\Omega_{\i}=[-4,\;4]\times[-4,\;4]$. We choose the width of the boundary domain $\Gamma$ to be 3.5.
In order to  examine  the absorbing property, we tracked the change of the total electrons in the subdomain $\Omega_{\i}$: \ie,  $N_I(t)=\int_{\Omega_\i} \delta n(\bd r, t) d\bd r$. As shown in Figure 
\ref{eq: ech}, the first-order BC offers reasonable accuracy at the initial period $0\le t \le 6$. However, the electrons will return to the subdomain later, and $N_I$ will rise significantly. 
The second-order ABCs yield much more accurate results.  Among the two second-order methods, the ABC IIb offers slightly better accuracy. In the construction of the ABCs,
we choose $s_0=0.2.$

  \begin{figure}[htbp]
\begin{center}
\includegraphics[width=2.6cm,height=16cm,angle=-90]{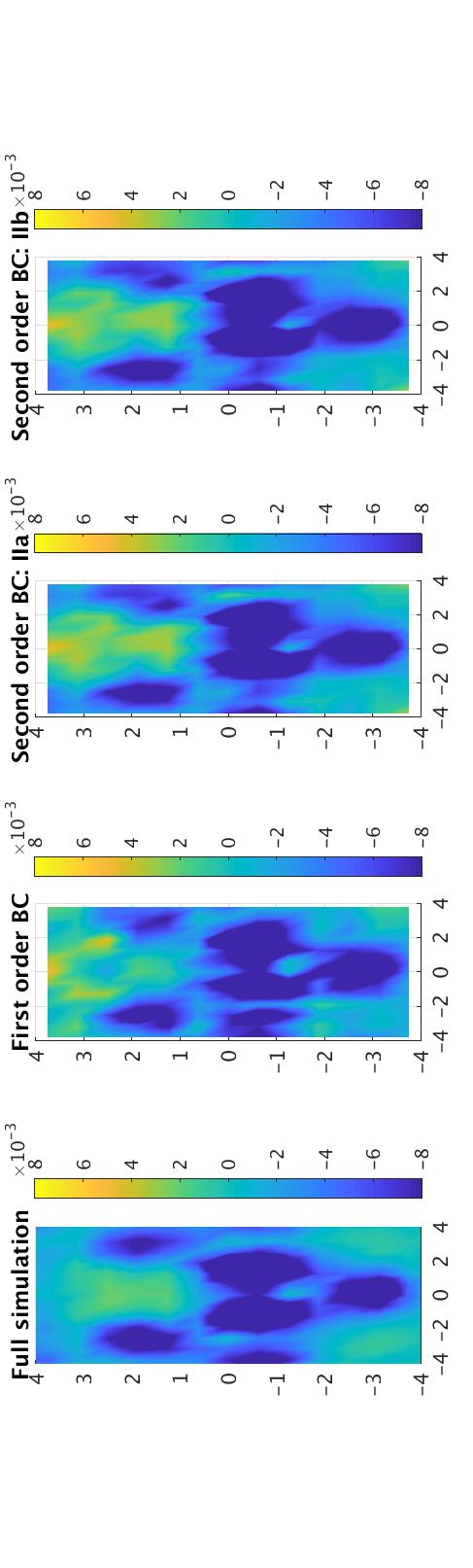}\\
\includegraphics[width=2.6cm,height=16cm,angle=-90]{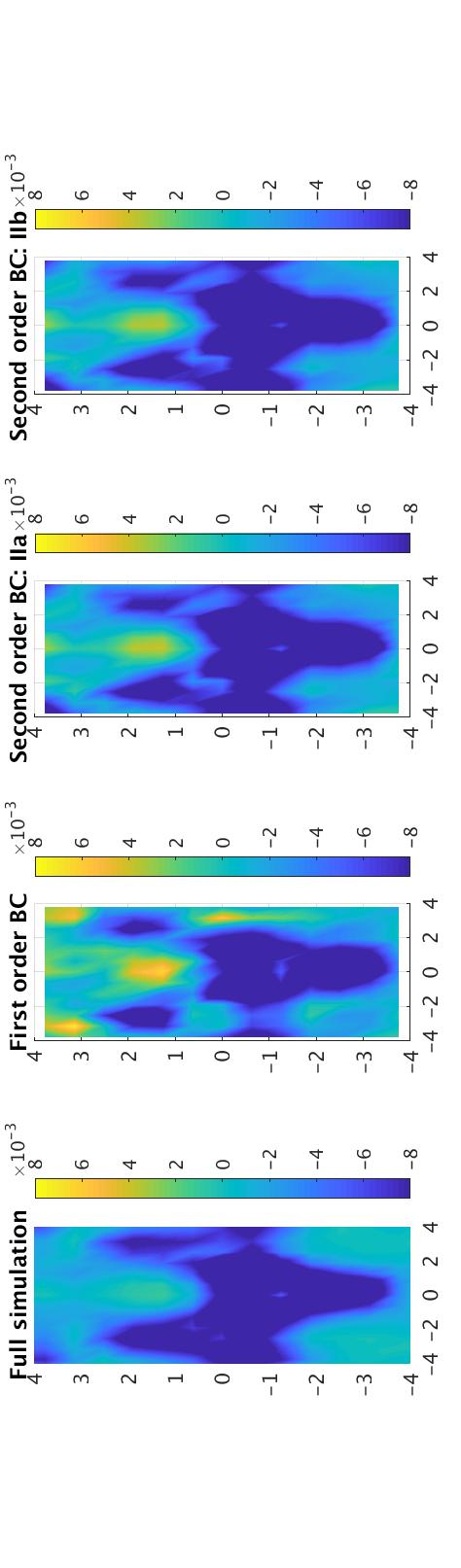}
\includegraphics[width=2.6cm,height=16cm,angle=-90]{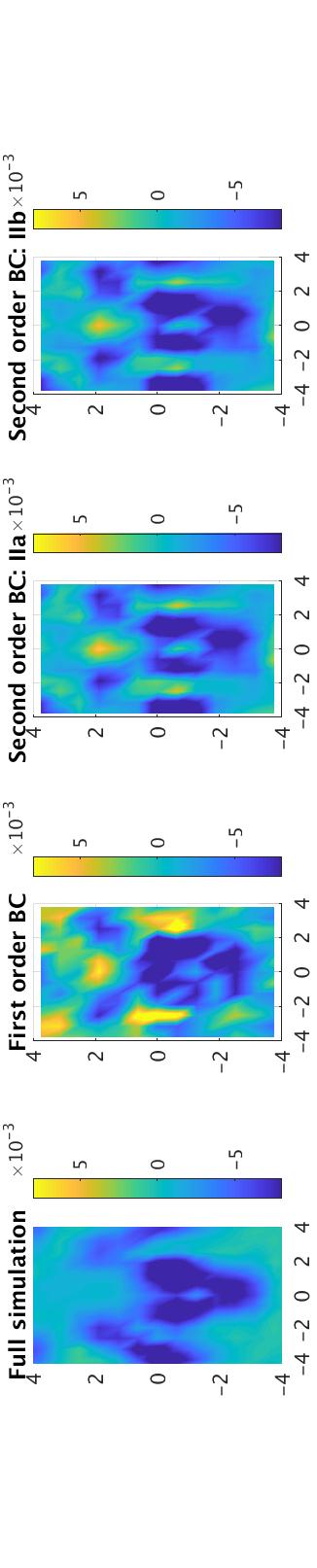}
\includegraphics[width=2.6cm,height=16cm,angle=-90]{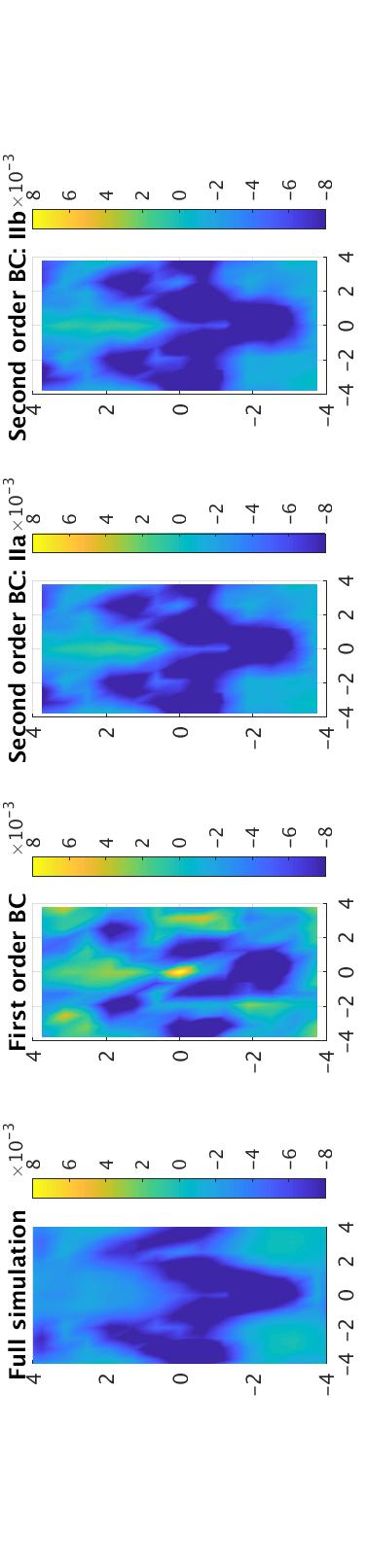}
\end{center}
\caption{Comparison of the ABCs.  From top to bottom: time = 6.4, time= 7.6, time = 8.8, time=10. }
\label{fig: abc}
\end{figure}

Finally, we compare the electron density distribution in the domain $\Omega_\i$,  at various points in time. The results are shown in Figure \ref{fig: abc}. 
We can observe that for the first-order ABC, there are significant reflections of the electron density (yellow peaks), while the second order ABCs are quite effective
in preventing such reflections.

\section{Summary and Discussions}\label{sum}
We have presented some absorbing boundary conditions for computer simulations based on models of time-dependent density-functional theory. The boundary conditions have been formulated using the density-matrix, and they are expressed in terms of the elements of the density-matrix at the boundary of the computational domain. A further reduction is introduced to approximate the memory term in the integral, which subsequently leads to modified density-matrix equations at the boundary. Numerical experiments have shown that these boundary conditions are very effective in suppressing boundary reflections. We would like to  further comment that while this procedure does lead to a modification of the Hamiltonian with a non-zero imaginary part at the boundary, it is more general and systematic than the  imaginary potential idea \cite{child1991analysis,muga2004complex}, in that the modified Hamiltonian is determined from a full model, and both the real and imaginary parts 
are different from the original matrix.  In fact, the analysis of classical absorbing boundary conditions suggests that it is often not sufficient to simply introduce damping. Rather, a perfect matching is only accomplished when  both the damping coefficient and the wave speed are modified, so that the impedance matches at the boundary \cite{berenger1994}. Our approach is parallel to the Dirchlet-to-Neumann map approach,
 \cite{givoli1998discrete} and it constructs the modified equations by directly eliminating the extra degrees of freedom in the exterior.

A natural extension of the current framework is to the linear response formalism \cite{casida1995time,muta2002solving,yabana2006real}, including both the Sternheimer and Casida's formulation for computing transport properties and excitation states. For those cases, the absorbing boundary conditions will be expressed in the frequency domain. Another important extension is to quantum molecular dynamics, where the nuclei are also allowed to move. In this case,  absorbing boundary conditions can also be designed to propagate out phonons that have been generated 
in the computational domain \cite{Wu2018}.  It would be interesting to formulate these two absorbing boundary conditions under the same framework. This work is underway.

\begin{acknowledgements}
 This research was supported by NSF under grant DMS-1522617 and DMS-1619661.
\end{acknowledgements}

\bibliographystyle{plain}

%


\end{document}